\input amstex
\input epsf
\documentstyle{amsppt}

\def\ls{\leqslant}
\def\gs{\geqslant}

\def\F{\Psi}
\def\f{\psi}
\def\S{S}
\TagsOnRight
\NoBlackBoxes

\topmatter
\title
Counting occurrences of $132$ in a permutation
\endtitle
\author Toufik Mansour$^*$ and Alek Vainshtein$^\dag$ \endauthor
\affil $^*$ Department of Mathematics\\
$^\dag$ Department of Mathematics and Department of Computer Science\\ 
University of Haifa, Haifa, Israel 31905\\ 
{\tt tmansur\@study.haifa.ac.il},
{\tt alek\@mathcs.haifa.ac.il}
\endaffil

\abstract
We study the
generating function for the number of permutations on $n$ letters containing
exactly $r\gs0$ occurrences of $132$. It is shown that finding this function
for a given $r$ amounts to a routine check of all permutations in $S_{2r}$.
\medskip
\noindent {\smc 2000 Mathematics Subject Classification}: 
Primary 05A05, 05A15; Secondary 05C90
\endabstract

\leftheadtext{Toufik Mansour and Alek Vainshtein}
\endtopmatter

\document
\heading 1. Introduction \endheading

Let $\pi\in \S_n$ and $\tau\in \S_{m}$ be two permutations. An 
{\it occurrence\/}
of $\tau$ in $\pi$ is a subsequence $1\ls i_1<i_2<\dots<i_{m}\ls n$ 
such that
$(\pi(i_1),\dots,\pi(i_{m}))$ is order-isomorphic to $\tau$; in such a 
context, $\tau$ is usually called a {\it pattern}. 

Recently, much attention has been paid to the problem of counting the number 
$\f_r^\tau(n)$ of
permutations of length $n$ containing a given number $r\gs0$ of occurrences
of a certain pattern $\tau$. Most of the authors consider only the case
$r=0$, thus studying permutations {\it avoiding\/} a given pattern. Only a 
few papers consider the case $r>0$, usually restricting themselves to the
patterns of length $3$. In fact, simple algebraic considerations show that 
there are only two essentially different cases for $\tau\in \S_3$, namely, 
$\tau=123$ and $\tau=132$. Noonan \cite{No} has proved that 
$\f_1^{123}(n)=\frac 3n\binom{2n}{n-3}$. A general approach to the problem was 
suggested by Noonan and Zeilberger \cite{NZ}; they gave another proof of
Noonan's result, and conjectured that 
$$
\f_2^{123}(n)=\frac{59n^2+117n+100}{2n(2n-1)(n+5)}\binom{2n}{n-4}
$$ 
and $\f_1^{132}(n)=\binom{2n-3}{n-3}$. The latter conjecture was proved by
B\'ona in \cite{B2}. A general conjecture of Noonan and Zeilberger states
that $\f_r^\tau(n)$ is $P$-recursive in $n$ for any $r$ and $\tau$. 
It was proved by B\'ona \cite{B1} for $\tau=132$. However, as stated 
in \cite{B1}, a challenging question is to describe $\f_r^{\tau}(n)$, 
$\tau\in \S_3$, explicitly for any given $r$. 

In this note we suggest a new approach to this problem in the case 
$\tau=132$, which allows to get
an explicit expression for $\f_r(n)=\f_r^{132}(n)$ for any given $r$. More 
precisely,
we present an algorithm that computes the generating function
$\F_r(x)=\sum_{n\gs0}\f_r(n)x^n$ for any $r\gs0$. To get the result for a
given $r$, the algorithm performs certain routine checks for each element
of the symmetric group $\S_{2r}$. The algorithm has been implemented in
C, and yielded explicit results for $1\ls r\ls 6$.

The authors are sincerely grateful to M.~Fulmek and A.~Robertson for inspiring
discussions.

\heading 2. Preliminary results \endheading

To any $\pi\in \S_n$ we assign a bipartite graph $G_\pi$ in the following way.
The vertices in one part of $G_\pi$, denoted $V_1$, are the entries of $\pi$,
and the vertices of the second part, denoted $V_3$, are the occurrences of
$132$ in $\pi$. Entry $i\in V_1$ is connected by an edge to occurrence 
$j\in V_3$ if $i$ enters $j$. For example, let $\pi=57614283$, then $\pi$
contains $5$ occurrences of $132$, and the graph $G_\pi$ is presented on
Figure~1.

\vskip 10pt
\centerline{\hbox{\epsfxsize=8cm\epsfbox{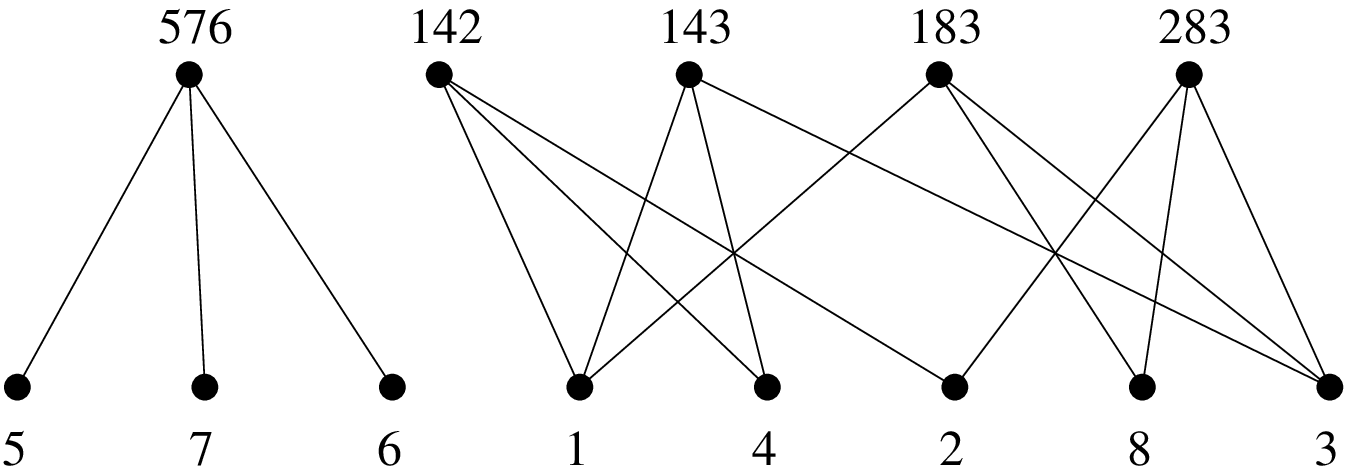}}}
\midspace{1mm}
\caption{Fig.~1. Graph $G_\pi$ for $\pi=57614283$}
\vskip 5pt

Let $\widetilde G$ be an arbitrary connected component of $G_\pi$, and let
$\widetilde V$ be its vertex set. We denote 
$\widetilde V_1=\widetilde V\cap V_1$, $\widetilde V_3=\widetilde V\cap V_3$,
$t_1=|\widetilde V_1|$, $t_3=|\widetilde V_3|$.

\proclaim{Lemma 1} For any connected component $\widetilde G$ of $G_\pi$ one
has $t_1\ls 2t_3+1$.
\endproclaim

\demo{Proof} Assume to the contrary that the above statement is not true.
Consider the smallest $n$ for which there exists $\pi\in \S_n$ such that
for some connected component  $\widetilde G$ of $G_\pi$ one has 
$$
t_1>2t_3+1. \tag $\ast$
$$
Evidently, $\widetilde G$ contains more than one vertex, since otherwise 
$t_1=1$, $t_3=0$, which contradicts  $(\ast)$. Let $l$ be the number of 
leaves in $\widetilde G$ (recall that a leaf is a vertex of
degree~$1$). Clearly, all the leaves belong to $\widetilde V_1$;
the degree of any other vertex in $\widetilde V_1$ is at least $2$, while
the degree of any vertex in $\widetilde V_3$ equals $3$. Calculating the 
number of edges in $\widetilde G$ by two different ways, we get 
$l+2(t_1-l)\ls 3t_3$, which together with ($\ast$) gives $l>t_3+2$.
This means that there exist two leaves $u,v\in \widetilde V_1$ incident to
the same vertex $a\in \widetilde V_3$. 

Let $w\in \widetilde V_1$ be the third vertex incident to $a$. If $w$ is a 
leaf, then $\widetilde G$ contains only four vertices $a,u,v,w$, and hence
$t_1=3$, $t_3=1$, which contradicts  $(\ast)$. Hence, the degree of $w$
is at least $2$. Delete the entries $u,v$ from $\pi$ and consider the
corresponding permutation $\pi'\in \S_{n-2}$. Denote by $\widetilde G'$ the
connected component of $G_{\pi'}$ containing $w$. Since the degree of $w$
in $\widetilde G$ was at least $2$, we see that $\widetilde G'$ is obtained
from $\widetilde G$ by deleting vertices $u$, $v$, and $a$. Therefore,
$t_1'=t_1-2$, $t_3'=t_3-1$, and hence $t_1'>2t_3'+1$,  a contradiction
to the minimality of $n$.
\qed
\enddemo

Denote by $G_\pi^n$ the connected component of $G_\pi$ containing entry $n$.
Let $\pi(i_1),\dots,\allowmathbreak\pi(i_s)$ be the entries of $\pi$ 
belonging to 
$G_\pi^n$, and let $\sigma=\sigma_\pi\in \S_s$ be the corresponding 
permutation.
We say that $\pi(i_1),\dots,\pi(i_s)$ is the {\it kernel\/} of $\pi$
and denote it $\ker\pi$; 
$\sigma$ is called the {\it shape\/} of the kernel, or the {\it kernel shape},
$s$ is called the 
{\it size\/} of the kernel, and the number of occurrences of $132$ in 
$\ker\pi$ is called the {\it capacity\/} of the kernel. For 
example, for $\pi=57614283$ as above, the kernel equals $14283$, its shape
is $14253$, the size equals $5$, and the capacity equals $4$.

The following statement is implied immediately by Lemma~1.

\proclaim{Theorem~1} Let $\pi\in \S_n$ contain exactly $r$ occurrences of $132$,
then the size of the kernel of $\pi$ is at most $2r+1$.
\endproclaim

We say that $\rho$ is a {\it kernel permutation\/} if it is the kernel shape
for some permutation $\pi$. Evidently $\rho$ is a kernel 
permutation if and only if $\sigma_\rho=\rho$.

Let $\rho\in \S_s$ be an arbitrary kernel permutation. We denote by 
$\S(\rho)$ the
set of all the permutations of all possible sizes whose kernel shape equals
$\rho$. For any $\pi\in \S(\rho)$ we define the {\it kernel cell 
decomposition\/} as follows. The number of cells in the decomposition equals
$s(s+1)$. Let $\ker\pi=\pi(i_1),\dots,\pi(i_s)$;
the {\it cell\/} $C_{ml}=C_{ml}(\pi)$ for $1\ls l\ls s+1$ and $1\ls m\ls s$
is defined by
$$
C_{ml}(\pi)=\{\pi(j)\: i_{l-1}<j<i_{l}, \;\pi(i_{\rho^{-1}(m-1)})<\pi(j)<
\pi(i_{\rho^{-1}(m)})\},
$$
where $i_0=0$, $i_{s+1}=n+1$,  and $\alpha(0)=0$ for any $\alpha$. If
$\pi$ coincides with 
$\rho$ itself, then all the cells in the decomposition are empty.
An arbitrary permutation in $\S(\rho)$ is obtained  by filling in some of
the cells in the cell decomposition. A cell $C$ is called 
{\it infeasible\/} if the existence of an entry $a\in C$ would imply
an occurrence of $132$ that contains $a$ and two other entries $x,y\in\ker\pi$.
Clearly, all infeasible cells are empty for any $\pi\in \S(\rho)$.
All the remaining cells are called {\it feasible\/}; a feasible cell may,
or may not, be empty. Consider the permutation $\pi=67382451$. The kernel
of $\pi$ equals $3845$, its shape is $1423$. The cell 
decomposition of $\pi$ contains four feasible cells: $C_{13}=\{2\}$,
$C_{14}=\varnothing$, $C_{15}=\{1\}$, and $C_{41}=\{6,7\}$, see Figure~2. 
All the other 
cells are infeasible; for example, $C_{32}$ is infeasible, since if 
$a\in C_{32}$, then $a\pi'(i_2)\pi'(i_4)$ is an occurrence of $132$
for any $\pi'$ whose kernel is of shape $1423$. 

\vskip 10pt
\centerline{\hbox{\epsfxsize=8cm\epsfbox{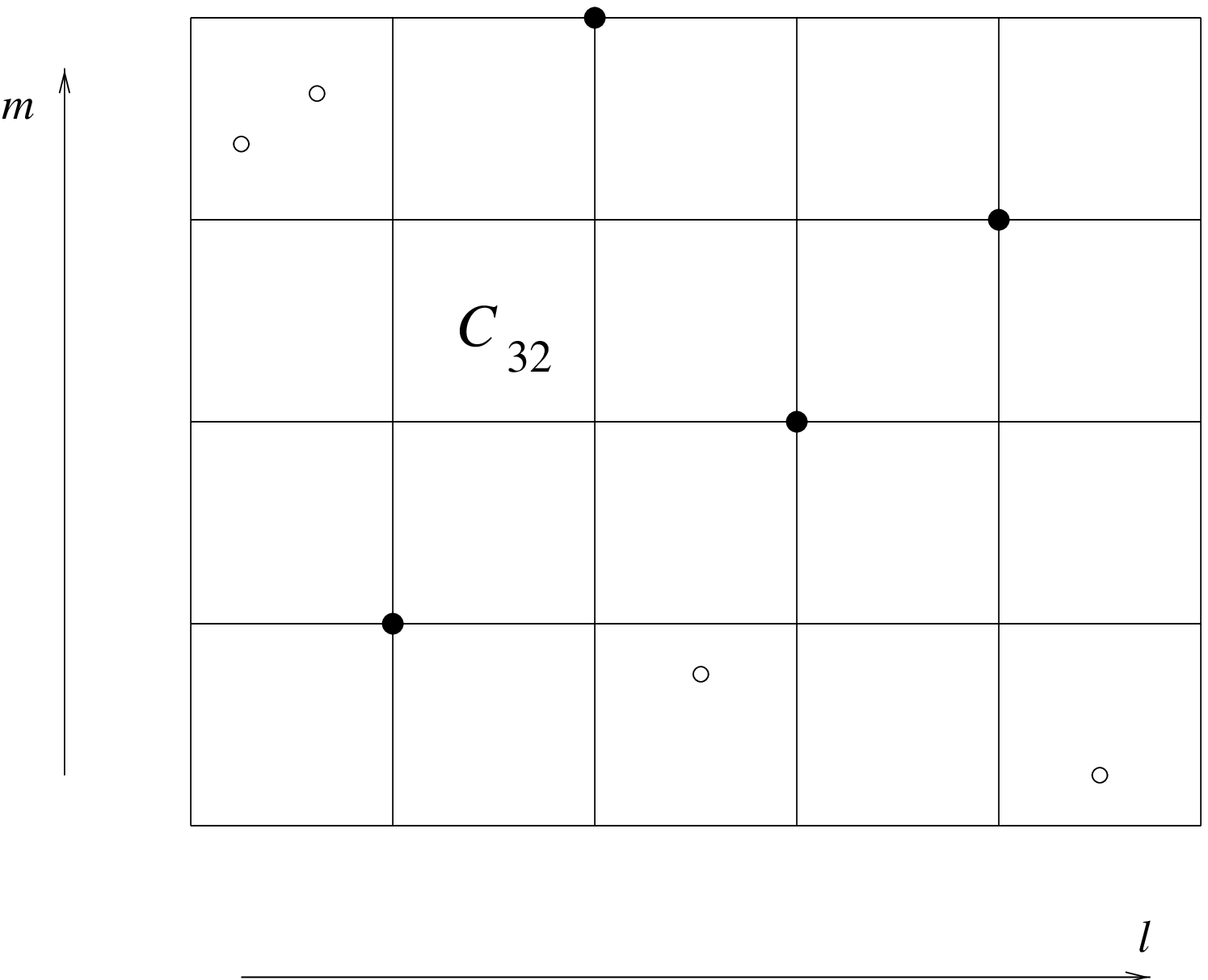}}}
\midspace{1mm}
\caption{Fig.~2. Kernel cell decomposition for $\pi\in\S(1423)$}
\vskip 5pt 

As another example, permutation
$\widetilde \pi=10\,11\,7\,12\,4\,6\,5\,8\,3\,9\,2\,1$ belongs to the same
class $\S(1423)$. Its kernel is $7\,12\,8\,9$, and the feasible cells are
$C_{13}=\{4,6,5\}$, $C_{14}=\{3\}$, $C_{15}=\{2,1\}$, $C_{41}=\{10,11\}$.

Given a cell $C_{ij}$ in the kernel cell decomposition, all the kernel
entries can be positioned with respect to $C_{ij}$. We say that
$x=\pi(i_k)\in\ker\pi$ lies {\it below\/} $C_{ij}$ if $\rho(k)<i$, and
{\it above\/} $C_{ij}$ if $\rho(k)\gs i$. Similarly, $x$ lies to the
{\it left\/} of $C_{ij}$ if $k<j$, and to the {\it right\/} of
$C_{ij}$ if $k \gs j$. As usual, we say that $x$ lies to the {\it
southwest\/} of $C_{ij}$ if it lies below $C_{ij}$ and to the left of
it; the other three directions, northwest, southeast, and northeast,
are defined similarly.

The following statement plays a crucial role in our considerations.

\proclaim{Lemma 2} Let $\pi\in\S(\rho)$ and
$\pi(i_k)\in\ker\pi$, then any cell $C_{ml}$ such that
$l> k$ and $m> \rho(k)$ is infeasible.
\endproclaim

\demo{Proof} Assume to the contrary that there exist $l$ and $m$ as above
such that $C_{ml}$ is feasible, and
consider the partition of the entries in $\ker\pi$
into two subsets: $\varkappa_1$ containing the entries of $\ker\pi$
that lie to the southwest of $C_{ml}$ 
and $\varkappa_2$ containing the rest of the entries. Observe that
$\pi(i_k)\ne n$, since $n\gs m>\rho(k)=\pi(i_k)$. Moreover, 
$\pi(i_k)\in\varkappa_1$ and $n\in\varkappa_2$. Since $\pi(i_k)$ and $n$ 
belong to the same 
connected component of $G_\pi$, there exists at least one occurrence of
$132$ whose elements are distributed between $\varkappa_1$ and $\varkappa_2$.
Let $a$ denote the minimal entry in this occurrence, let $c$ denote its maximal
entry, and let $b$ denote the remaining entry.

Evidently, $a\in\varkappa_1$. Assume first that $c\in\varkappa_2$.

If $c$ lies to the left of $C_{ml}$, then the existence of $z\in C_{ml}$
would imply that $acz$ is an occurrence of $132$, and hence $C_{ml}$ is 
infeasible.

If $c$ lies to the northeast of $C_{ml}$ and $b$ lies above $C_{ml}$,
then the existence of $z\in C_{ml}$
would imply that $zcb$ is an occurrence of $132$, and hence $C_{ml}$ is 
infeasible. If $b$ lies below $C_{ml}$,
then the existence of $z\in C_{ml}$
would imply that $azb$ is an occurrence of $132$, and hence $C_{ml}$ is 
infeasible.

If $c$ lies to the southeast of $C_{ml}$, then the existence of 
$z\in C_{ml}$
would imply that $azc$ is an occurrence of $132$, and hence $C_{ml}$ is 
infeasible.

It remains to consider the case $c\in\varkappa_1$, which means that
$b$ belongs to $\varkappa_2$ and lies to the southeast of $C_{ml}$.
Hence, the existence of $z\in C_{ml}$
would imply that $azb$ is an occurrence of $132$, and hence $C_{ml}$ is 
infeasible.
\qed
\enddemo

As an easy corollary of Lemma~2, we get the following proposition.
Let us define a partial order $\prec$ on the set of all feasible cells
by saying that $C_{ml}\prec C_{m'l'}\ne C_{ml}$ if $m\gs m'$ and $l\ls l'$.

\proclaim{Lemma 3} $\prec$ is a linear order.
\endproclaim

\demo{Proof} Assume to the contrary that there exist two feasible cells
$C_{ml}$ and $C_{m'l'}$ such that $l<l'$ and $m<m'$, and consider the entry
$x=\pi(i_{l})\in\ker\pi$. By Lemma~2, $x>\pi(i_{\rho^{-1}(m'-1)})$, that is, 
$x$ lies above the cell $C_{m'l'}$, since otherwise $C_{m'l'}$ would be
infeasible. For the same reason, $y=\pi(i_{\rho^{-1}(m'-1)})$ lies to the 
right of $C_{m'l'}$, and hence to the right of $x$. Therefore, the
existence of $z\in C_{ml}$ would imply that $zxy$ is an occurrence of $132$,
and hence $C_{ml}$ is infeasible, a contradiction.
\qed
\enddemo

Consider now the dependence between two nonempty feasible cells lying
on the same horizontal or vertical level.

\proclaim{Lemma 4} Let $C_{ml}$ and $C_{ml'}$ be two nonempty feasible cells
such that $l<l'$. Then for any pair of entries $a\in C_{ml}$, $b\in C_{ml'}$,
one has $a>b$.
\endproclaim

\demo{Proof} Assume to the contrary that there exists a pair $a\in C_{ml}$,
$b\in C_{ml'}$ such that $a<b$. Consider the entry $x=\pi(i_{l})\in\ker\pi$.
By Lemma~2, $x>b$, since otherwise $C_{ml'}$ would be infeasible. Hence 
$axb$ is an occurrence of $132$, which means that both $a$ and $b$ belong to
$\ker\pi$, a contradiction.
\qed
\enddemo

\proclaim{Lemma 5}  Let $C_{ml}$ and $C_{m'l}$ be two nonempty feasible cells
such that $m<m'$. Then any  entry $a\in C_{ml}$ lies to the right of any
entry $b\in C_{m'l}$.
\endproclaim

\demo{Proof} Assume to the contrary that there exists a pair $a\in C_{ml}$,
$b\in C_{m'l}$ such that $a$ lies to the left of $b$. Consider the entry 
$y=\pi(i_{\rho^{-1}(m'-1)})\in\ker\pi$.
By Lemma~2, $y$ lies to the right of $b$, since otherwise $C_{m'l}$ would be 
infeasible. Hence $aby$ is an occurrence of $132$, which means that both $a$ 
and $b$ belong to $\ker\pi$, a contradiction.
\qed
\enddemo

Lemmas~3--5 yield immediately the following two results.

\proclaim{Theorem 2} Let $\widetilde G$ be a connected component of
$G_\pi$ distinct from $G_\pi^n$. Then all the vertices in $\widetilde V_1$
belong to the same feasible cell in the kernel cell decomposition of $\pi$.
\endproclaim

Let $F(\rho)$ be the set of all feasible cells in the kernel cell decomposition
corresponding to permutations in $\S(\rho)$, and let $f(\rho)=|F(\rho)|$.
We denote the cells in $F(\rho)$ by $C^1,\dots,C^{f(\rho)}$ in such a way
that $C^i\prec C^j$ whenever $i<j$.

\proclaim{Theorem 3} For any given sequence  $\alpha_1,\dots,\alpha_{f(\rho)}$
of arbitrary permutations there exists $\pi\in \S(\rho)$ such that the content
of $C^i$ is order-isomorphic to $\alpha_i$.
\endproclaim
 
\heading 3. Main Theorem  and explicit results \endheading

Let $\rho$ be a kernel permutation, and let $s(\rho)$, $c(\rho)$, and 
$f(\rho)$ be the size of $\rho$, the capacity of $\rho$, and the number of
feasible cells in the cell decomposition associated with $\rho$, respectively.
Denote by $K$ the set of all kernel permutations, and by $K_t$ the set of all
kernel shapes for  permutations in $\S_t$. The main result of 
this note
can be formulated as follows.

\proclaim{Theorem 4} For any $r\gs1$,
$$
\F_r(x)=\sum_{\rho\in K_{2r+1}}\left(x^{s(\rho)}\sum_{r_1+\dots+r_{f(\rho)}=r-c(\rho)}
\;\prod_{j=1}^{f(\rho)}\F_{r_j}(x)\right), \tag $\ast\ast$
$$
where $r_j\gs0$ for $1\ls j\ls f(\rho)$.
\endproclaim

\demo{Proof} For any $\rho\in K$, denote by $\F_r^\rho(x)$ the generating
function for the number of permutations in $\pi\in\S_n\cap \S(\rho)$ containing
exactly $r$ occurrences of $132$. Evidently, $\F_r(x)=\sum_{\rho\in K}
\F_r^\rho(x)$. To find $\F_r^\rho(x)$, recall that the kernel of any
$\pi$ as above contains exactly $c(\rho)$ occurrences of $132$. The remaining
$r-c(\rho)$ occurrences of $132$ are distributed between the feasible cells 
of the kernel cell decomposition of $\pi$. By Theorem~2, each occurrence
of $132$ belongs entirely to one feasible cell. Besides, it follows from
Theorem~3, that occurrences of $132$ in different cells do not influence
one another. Therefore,
$$
\F_r^\rho(x)=x^{s(\rho)}\sum_{r_1+\dots+r_{f(\rho)}=r-c(\rho)}
\;\prod_{j=1}^{f(\rho)}\F_{r_j}(x),
$$
and we get the expression similar to ($\ast\ast$) with the only difference
that the outer sum is over all $\rho\in K$. However, if $\rho\in K_t$ for
$t>2r+1$, then by Theorem~1, $c(\rho)>r$, and hence $\F_r^\rho(x)\equiv 0$.
\qed
\enddemo

Theorem~4 provides a finite algorithm for finding $\F_r(x)$ for any given
$r>0$, since we have to consider all permutations in $\S_{2r+1}$, and to 
perform certain routine operations with all shapes found so far. Moreover,
the amount of search can be decreased substantially due to the following 
proposition.

\proclaim{Proposition} The only kernel permutation of capacity $r\gs1$ and
size $2r+1$ is $2r-1\,2r+1\,2r-3\,2r\,\dots 2r-2j-3\, 2r-2j\,\dots 1\,4\,2$.
Its contribution to $\F_r(x)$ equals $x^{2r+1}\F_0^{r+2}(x)$.
\endproclaim

This proposition is proved easily by induction, similarly to
Lemma~1. The feasible cells in the corresponding cell decomposition
are $C_{2r-2j+1,2j+1}$, $j=0,\dots,r$, and $C_{1,2r+2}$, hence the
contribution to $\F_r(x)$ is as described.

By the above proposition, it suffices to search only permutations in $\S_{2r}$.
Below we present several explicit calculations.

Let us start from the case $r=0$. Observe that ($\ast\ast$) remains valid for 
$r=0$, provided the left hand side is replaced by $\F_r(x)-1$; subtracting
$1$ here accounts for the empty permutation. So, we begin with finding 
kernel shapes for all permutations in $\S_1$. The only shape obtained is 
$\rho_1=1$, and it is easy to see that $s(\rho_1)=1$, $c(\rho_1)=0$, and 
$f(\rho_1)=2$
(since both cells $C_{11}$ and $C_{12}$ are feasible). Therefore, we get
$\F_0(x)-1=x\F_0^2(x)$, which means that
$$
\F_0(x)=\frac{1-\sqrt{1-4x}}{2x},
$$
the generating function of Catalan numbers. 

Let now $r=1$. Since permutations in $\S_2$ do not exhibit kernel shapes 
distinct from $\rho_1$, the only possible new shape is the exceptional
one, $\rho_2=132$, whose contribution equals $x^3\F_0^3(x)$. 
Therefore, ($**$) amounts to
$$
\F_1(x)=2x\F_0(x)\F_1(x) +x^3\F_0^3(x),
$$
and we get the following result.

\proclaim{Corollary 1 ({\rm B\'ona \cite{B2, Theorem~5}})}
$$
\F_1(x)=\frac12\left(x-1+(1-3x)(1-4x)^{-1/2}\right);
$$
equivalently,
$$
\f_1(n)=\binom{2n-3}{n-3}
$$
for $n\gs3$.
\endproclaim

Let $r=2$. We have to check the kernel shapes of permutations in $\S_4$.
Exhaustive search adds four new shapes to the previous list; these are
$1243$, $1342$, $1423$, and $2143$; besides, there is the exceptional 
$35142\in\S_5$. Calculation of the parameters
$s$, $c$, $f$ is straightforward, and we get

\proclaim{Corollary 2}
$$
\F_2(x)=\frac12\left(x^2+3x-2+(2x^4-4x^3+29x^2-15x+2)(1-4x)^{-3/2}\right);
$$
equivalently,
$$
\f_2(n)=\frac{n^3+17n^2-80n+80}{2n(n-1)}\binom{2n-6}{n-2}
$$
for $n\gs4$.
\endproclaim 

Let $r=3,4,5,6$; exhaustive search in $\S_6$, $\S_8$, $\S_{10}$, and $\S_{12}$ 
reveals $20$, $104$, $503$, and $2576$ new nonexceptional kernel shapes, 
respectively, and we get

\proclaim{Corollary 3} Let $3\ls r\ls 6$, then
$$
\F_r(x)=\frac12\left(P_r(x)+Q_r(x)(1-4x)^{-r+1/2}\right),
$$
where
$$\align
P_3(x)&=2x^3-5x^2+7x-2,\\
P_4(x)&=5x^4-7x^3+2x^2+8x-3,\\
P_5(x)&=14x^5-17x^4+x^3-16x^2+14x-2\endalign
$$
and
$$\align
Q_3(x)&=-22x^6-106x^5+292x^4-302x^3+135x^2-27x+2,\\
Q_4(x)&=2x^9+218x^8+1074x^7-1754x^6+388x^5+1087x^4\\
&\qquad\qquad\qquad\qquad-945x^3+320x^2-50x+3,\\
Q_5(x)&=-50x^{11}-2568x^{10}-10826x^9+16252x^8-12466x^7+16184x^6-16480x^5\\
&\qquad\qquad\qquad\qquad +9191x^4-2893x^3+520x^2-50x+2.\endalign
$$
Equivalently,
$$
\f_r(n)=R_r(n)\frac{(2n-3r)!}{n!r!(n-r-2)!},
$$
for $n\gs r+2$, where
$$\align
R_3(n)&=n^6+51n^5-407n^4-99n^3+7750n^2-22416n+20160,\\
R_4(n)&=n^9+102n^8-282n^7-12264n^6+32589n^5+891978n^4\\
&\qquad\qquad -7589428n^3+25452024n^2-39821760n+23950080,\\
R_5(n)&=n^{12}+170n^{11}+1861n^{10}-88090n^9-307617n^8+27882510n^7\\
&\qquad\qquad
-348117457n^6+2119611370n^5-6970280884n^4\\
&\qquad\qquad
+10530947320n^3+2614396896n^2
-30327454080n+29059430400.
\endalign
$$

The expressions for $P_6(x)$, $Q_6(x)$, and $R_6(n)$ are too long to be 
presented here.
\endproclaim 

\heading 4. Further results and open questions \endheading

As an easy consequence of Theorem~4 we get the following results due
to B\'ona \cite{B1}.

\proclaim{Corollary 4} Let $r\gs0$, then $\F_r(x)$ is a rational function in
the variables $x$ and $\sqrt{1-4x}$. 
\endproclaim

In fact, B\'ona has proved a stronger result, claiming that
$$
\F_r(x)=P_r(x)+Q_r(x)(1-4x)^{-r+1/2}, \tag $***$
$$
where $P_r(x)$ and $Q_r(x)$ are polynomials and
$1-4x$ does not divide $Q_r(x)$. We were unable to prove this result; however,
it stems almost immediately from the following conjecture.

\proclaim{Conjecture 1} For any kernel permutation $\rho\ne1$,
$$
s(\rho)\gs f(\rho).
$$
\endproclaim

Indeed, it is easy to see that $\F_r(x)$ enters the right hand side
of ($**$) with the coefficient $2x\F_0(x)$, which is a partial contribution 
of the
kernel shape $\rho_1=1$. Since $1-2x\F_0(x)=\sqrt{1-4x}$, we get 
by induction from
($**$) that $\sqrt{1-4x}\F_r(x)$ equals the sum of fractions whose 
denominators are of the form $x^{d}(1-4x)^{r-c(\rho)-f(\rho)/2}$, where 
$d\ls f(\rho)$. On the other hand, each fraction is multiplied by
$x^{s(\rho)}$, hence if $s(\rho)\gs f(\rho)$ as conjectured, then 
$x^d$ in the denominator is cancelled. 
The maximal degree of $(1-4x)$ is attained for $\rho=\rho_1$,
and is equal to $r-1$, and we thus arrive at ($***$).

In view of our explicit results, we have even a stronger conjecture.

\proclaim{Conjecture 2} The polynomials $P_r(x)$ and $Q_r(x)$ in $(***)$
have halfinteger coefficients.
\endproclaim

Another direction would be to match the approach of this note with the
previous results on restricted $132$-avoiding permutations. Let $\Phi_r(x;k)$
be the generating function for the number of permutations in $\S_n$ containing
$r$ occurrences of $132$ and avoiding the pattern $12\dots k\in\S_k$. It was 
shown previously that $\Phi_r(x;k)$ can be expressed via Chebyshev polynomials
of the second kind for $r=0$ (\cite{CW}) and $r=1$ (\cite{MV}). Our new 
approach allows to get a recursion for $\Phi_r(x;k)$ for any given
$r\gs0$.

Let $\rho$ be a kernel permutation, and assume that the feasible cells of
the kernel cell decomposition associated with $\rho$ are ordered linearly
according to $\prec$. We denote by $l_j(\rho)$ the length of the longest
increasing subsequence of $\rho$ that lies to the north-east from $C^j$. 
For example, let $\rho=1423$, as on Figure~2. Then $l_1(\rho)=1$, 
$l_2(\rho)=2$, $l_3(\rho)=1$, $l_4(\rho)=0$.

\proclaim{Theorem 5} For any $r\gs1$ and $k\gs3$, 
$$
\Phi_r(x;k)=\sum_{\rho\in K_{2r+1}}\left(x^{s(\rho)}
\sum_{r_1+\dots+r_{f(\rho)}=r-c(\rho)}
\;\prod_{j=1}^{f(\rho)}\Phi_{r_j}(x;k-l_j(\rho))\right), 
$$
where $r_j\gs0$ for $1\ls j\ls f(\rho)$ and $\Phi_r(x;m)\equiv0$ for $m\ls0$.
\endproclaim

As in the case of $\F_r(x)$, the statement of the theorem remains valid for 
$r=0$, provided the left hand side is replaced by $\Phi_r(x;k)-1$. This allows
to recover known explicit expressions for $\Phi_r(x;k)$ for $r=0,1$, and to get
an expression for $r=2$, which is too long to be presented here.

This approach can be extended even further, to cover also permutations
containing $r$ occurrences of $132$ and avoiding other permutations in $\S_k$,
for example, $23\dots k1$.


\Refs
\widestnumber\key{CW}

\ref \key B1        
\by M.~Bona
\paper The number of permutations with exactly $r$ $132$-subsequences is
$P$-recursive in the size!
\jour Adv. Appl. Math.
\vol 18 \yr 1997 \pages 510--522
\endref

\ref \key B2
\by M.~Bona
\paper Permutations with one or two $132$-subsequences
\jour Discr. Math
\yr 1998 \vol 181 \pages 267--274
\endref

\ref \key CW
\by T.~Chow and J.~West
\paper Forbidden subsequences and Chebyshev polynomials
\jour Discr. Math.
\vol 204 \yr 1999 \pages 119--128
\endref

\ref\key MV
\by T.~Mansour and A.~Vainshtein
\paper Restricted permutations, continued fractions, and Chebyshev polynomials
\jour Electron. J. Combin.
\vol 7 \yr 2000 \finalinfo \#R17
\endref

\ref \key No
\by J.~Noonan
\paper The number of permutations containing exactly one increasing 
subsequence of length three
\jour Discr. Math.
\vol 152 \yr 1996 \pages 307--313
\endref

\ref \key NZ              
\by J.~Noonan and D.~Zeilberger
\paper The enumeration of permutations with a prescribed number of
``forbidden'' patterns
\jour Adv. Appl. Math.
\vol 17 \yr 1996 \pages 381--407
\endref


\endRefs
\enddocument